\documentclass[11pt,reqno]{amsart}
\usepackage{amsfonts,amssymb,amsmath}

\newtheorem{thm}{Theorem}

\newtheorem{cor}[thm]{Corollary}
\newtheorem{lem}[thm]{Lemma}

\DeclareMathOperator{\Ant}{\mathfrak{A}}
\DeclareMathOperator{\bmap}{\mathfrak b}

\numberwithin{equation}{section} \numberwithin{thm}{section}

\title{Maps on posets, and blockers}

\date{}
\author{Andrey O. Matveev}
\address{Data-Center Co., RU-620034, Ekaterinburg,
P.O.~Box~5, Russian~Federation} \email{aomatveev@\{dc.ru,
hotmail.com\}}

\begin{document}

\begin{abstract}
An order-theoretic generalization of Seymour relations describing
the connection between the set-theoretic blocker, deletion, and
contraction maps on clutters, is presented.
\end{abstract}

\maketitle

\pagestyle{myheadings} \markboth{Maps on posets, and
blockers}{Andrey~O.~Matveev}

\section{Introduction}
\thispagestyle{empty}

The aim of this note is to present an order-theoretic
generalization of Seymour relations~\cite{PDS} which describe the
set-theoretic blocker, deletion, and contraction maps on clutters,
see~(\ref{Label6}) below. Those relations are a powerful tool of
discrete mathematics, see, e.g.,~\cite{C,CH}.

A set $H$ is called a {\em blocking set\/} ({\em cover, system of
representatives, transversal}) for a nonempty family
$\mathcal{G}=\{G_1,\ldots,G_m\}$ of nonempty subsets of a finite
set if it holds $|H\cap G_k|>0$, for each $k\in\{1,\ldots,m\}$.
The family of all inclusion-minimal blocking sets for
$\mathcal{G}$ is called the {\em blocker\/} of $\mathcal{G}$, see,
e.g.,~\cite[Chapter~8]{GLS}. We denote the blocker of
$\mathcal{G}$ by $\mathcal{B}(\mathcal{G})$.

A family of subsets of a finite {\em ground set\/} $S$ is called a
{\em clutter} or a {\em Sperner family\/} if no set from that
family contains another. The empty clutter $\emptyset$ containing
no subsets of $S$, and the clutter $\{\hat{0}\}$ whose unique set
is the empty subset $\hat{0}$ of $S$, are called the {\em trivial
clutters} on $S$. The set-theoretic {\em blocker map\/} assigns to
a nontrivial clutter its blocker, and this map alternates the
trivial clutters: $\mathcal{B}(\emptyset):=\{\hat{0}\}$ and
$\mathcal{B}(\{\hat{0}\}):=\emptyset$, see, e.g.,~\cite{CFM}.

Let $X\subseteq S$. The set-theoretic {\em deletion $(\backslash
X)$, and contraction $(/X)$ maps} on clutters are defined in the
following way: if $\mathcal{G}$ is a nontrivial clutter on $S$
then the {\em deletion\/} $\mathcal{G}\backslash X$ is the family
$\{G\in\mathcal{G}:\ |G\cap X|=0\}$, and the {\em contraction\/}
$\mathcal{G}/X$ is the family of all inclusion-minimal sets from
the family $\{G-X:\ G\in\mathcal{G}\}$. One often says that the
clutters $\mathcal{G}\backslash X$ and $\mathcal{G}/ X$ are those
on the ground set $S-X$. The trivial clutters do not change under
the {\em deletion\/} and {\em contraction maps}:
$\emptyset\backslash X=\emptyset/X:=\emptyset$ and
$\{\hat{0}\}\backslash X=\{\hat{0}\}/X:=\{\hat{0}\}$.

Let $\mathcal{G}$ be a clutter on the ground set $S$.  We have
\begin{equation}
\label{Label5} \mathcal{B}(\mathcal{B}(\mathcal{G}))=\mathcal{G}\
,
\end{equation}
see~\cite{EF,L}; given a subset $X\subseteq S$, it holds
(see~\cite{PDS}):
\begin{equation} \label{Label6} \mathcal{B}(\mathcal{G})\backslash
X=\mathcal{B}(\mathcal{G}/X) \text{\ \ and\ \ }
\mathcal{B}(\mathcal{G})/ X=\mathcal{B}(\mathcal{G}\backslash X)\
.
\end{equation}

\section{A generalization of relations~(\ref{Label6})}

We refer the reader to~\cite[Chapter~3]{St} for information and
terminology in the theory of posets. See,
e.g.,~\cite[Chapter~IV]{Aigner} on the Galois correspondence and
(co)closure operators.

\begin{thm}
\label{Label7} Let $L$ be a finite poset. Let $\delta:L\to L$ be
an order-preserving map, and let $\gamma:L\to L$ be an
order-preserving map such that
\begin{equation}
\label{Label1} \gamma(x)\ \geq\ x\ ,
\end{equation}
for all $x\in L$. Let $\beta:L\to L$ be an order-reversing map
such that
\begin{equation}
\label{Label3} \beta(\beta(x))\ \geq\ x\ ,
\end{equation}
for all $x\in L$. Either of the relations {\rm(}for all $x\in L${\rm):}
\begin{equation}
\label{Label9} \beta(\delta(\beta(x)))\ \geq\ \gamma(x)\ ,
\end{equation}
\begin{equation}
\label{Label11} \beta(\gamma(\beta(x)))\ \geq\ \delta(x)\
\end{equation}
implies
\begin{equation}
\label{Label8} \delta(\beta(z))\ \leq\ \beta(\gamma(z))\ \leq\
\beta(z)\ \leq\ \gamma(\beta(z))\ \leq\ \beta(\delta(z))\ ,
\end{equation}
for any $z\in L$. Moreover, if $\beta(\beta(x))=x$, for all $x\in
L$, then either of the equalities
$\beta(\delta(\beta(x)))=\gamma(x)$ and
$\beta(\gamma(\beta(x)))=\delta(x)$, for all $x\in L$, implies
\begin{equation}\label{Label12}
\delta(\beta(z))\ =\ \beta(\gamma(z))\ \leq\ \beta(z)\ \leq\
\gamma(\beta(z))\ =\ \beta(\delta(z))\ ,
\end{equation}
for any $z\in L$.
\end{thm}

\begin{proof}
Relation~(\ref{Label1}) implies
\begin{equation*}
\beta(\gamma(z))\ \leq\ \beta(z)
\end{equation*}
because the map $\beta$ is order-reversing; moreover, we have
\begin{equation*}
\beta(z)\ \leq\ \gamma(\beta(z))\ .
\end{equation*}

We now prove
implication~(\ref{Label9})$\Longrightarrow$(\ref{Label8}).

On the one hand, with respect to~(\ref{Label3}), we have
$\delta(\beta(z))\leq\beta (\beta(\delta(\beta(z))))$. On the
other hand, since $\beta$ is order-reversing,
relation~(\ref{Label9}) implies\\ $\beta
(\beta(\delta(\beta(z))))\leq\beta (\gamma(z))$. We obtain
\begin{equation}
\label{Label18} \delta(\beta(z))\ \leq\ \beta(\gamma(z))\ .
\end{equation}

Further, on the one hand, relation~(\ref{Label9}) implies
$\gamma(\beta(z))\leq \beta(\delta(\beta(\beta(z))))$. On the
other hand, since $\beta(\beta(z))\geq z$ by~(\ref{Label3}), and
$\delta$ is order-preserving, and $\beta$ is order-reversing, we
obtain $\beta(\delta(\beta(\beta(z))))\leq \beta(\delta(z))$. We
conclude that
\begin{equation}
\label{Label19} \gamma(\beta(z))\ \leq\ \beta(\delta(z))\ ,
\end{equation}
and we are done.

We now prove
implication~(\ref{Label11})$\Longrightarrow$(\ref{Label8}).

On the one hand, with respect to~(\ref{Label11}), we have
$\delta(\beta(z))\leq\beta (\gamma(\beta(\beta(z))))$. On the
other hand, since $\beta$ is order-reversing, and $\gamma$ is
order-preserving, relation~(\ref{Label3}) implies $\beta
(\gamma(\beta(\beta(z))))\leq\beta (\gamma(z))$. We
obtain~(\ref{Label18}).

Further, on the one hand, relation~(\ref{Label3}) implies
$\gamma(\beta(z))\leq \beta(\beta(\gamma(\beta(z))))$. On the
other hand, since $\beta$ is order-reversing,
relation~(\ref{Label11}) implies\\
$\beta(\beta(\gamma(\beta(z))))\leq \beta(\delta(z))$. We come
to~(\ref{Label19}), and we are done.

The proof of relation~(\ref{Label12}) is now straightforward, with
respect to the argument above.
\end{proof}

Note that since the map $\beta$ in Theorem~\ref{Label7} is
order-reversing, and (\ref{Label3}) holds, it is a consequence of
\cite[Proposition~4.36(iii)]{Aigner} that we have
\begin{equation}
\label{Label10} \beta(\beta(\beta(x)))=\beta(x)\ ,
\end{equation}
for any $x\in L$.

To illustrate Theorem~\ref{Label7}, we give a comment
to~(\ref{Label6}). Let $P$ be a finite bounded poset of
cardinality greater than one, whose least element is denoted
$\hat{0}_P$. We denote by $\mathfrak{I}(A)$ and $\mathfrak{F}(A)$
the order ideal and filter of $P$ generated by an antichain
$A\subset P$, respectively. The {\em atoms\/} of $P$ are the
elements covering $\hat{0}_P$; we denote the set of all atoms of
$P$ by $P^{\mathrm a}$.

The antichains in $P$ compose a distributive lattice, denoted
$\Ant(P)$. In the present note, the antichains are
ordered in the following way: if $A',A''\in\Ant(P)$ then we set
\begin{equation*}
A'\leq A''\text{\ \ iff\ \ } \mathfrak{F}(A')\subseteq\mathfrak{F}(A'')\ .
\end{equation*}
We call the least element $\hat{0}_{\Ant(P)}$ and greatest element
$\hat{1}_{\Ant(P)}$ of $\Ant(P)$ the {\em trivial antichains\/} in
$P$ because, in the context of the present note, those antichains
are counterparts of the trivial clutters. Here $\hat{0}_{\Ant(P)}$
is the empty antichain in $P$, and $\hat{1}_{\Ant(P)}$ is the
one-element antichain $\{\hat{0}_P\}$.

\begin{itemize}
\item
If $\{a\}$ is a nontrivial one-element antichain in $P$ then the
order-theoretic {\em blocker $\bmap(a)$} of $\{a\}$ in $P$ is the
antichain
\begin{equation*}
\bmap(a):=\mathfrak{I}(a)\cap P^{\mathrm a}\ .
\end{equation*}

\item
If $A$ is a nontrivial antichain in $P$ then the order-theoretic
{\em blocker} $\bmap(A)$ of $A$ in $P$ is the following meet in
$\Ant(P)$:
\begin{equation}
\label{Label16} \bmap(A):=\bigwedge_{a\in A}\bmap(a)\ .
\end{equation}

\item
The order-theoretic {\em blockers\/} of the trivial antichains in
$P$ are:
\begin{equation*}
\bmap(\hat{0}_{\Ant(P)}):=\hat{1}_{\Ant(P)}\ ,\ \ \
\bmap(\hat{1}_{\Ant(P)}):=\hat{0}_{\Ant(P)}\ .
\end{equation*}
\end{itemize}
See~\cite{BjH,BPS,M1,M2,M3} on blockers in posets.

The map $\bmap:\Ant(P)\to\Ant(P)$ is called the order-theoretic
{\em blocker map\/} on $\Ant(P)$. That map is order-reversing,
with the property $\bmap(\bmap(A))\geq A$, for all $A\in\Ant(P)$.
Equality~(\ref{Label10}) implies
\begin{equation*}
\bmap(\bmap(\bmap(A)))=\bmap(A)\ ,
\end{equation*}
cf.~(\ref{Label5}). The posets with the property
$\bmap(\bmap(A))=A$, for all $A\in\Ant(P)$, are characterized
in~\cite{BjH}.


Let $X\subseteq P^{\mathrm a}$.

\begin{itemize}
\item
If $\{a\}$ is a nontrivial one-element antichain in $P$ then the
order-theoretic {\em deletion\/} $\{a\}\backslash X$ and {\em
contraction\/} $\{a\}/ X$ of $\{a\}$ in $P$ are the antichains
\begin{align*}
\{a\}\backslash X &:= \left\{\begin{array}{ll}\{a\}, & \mbox{if }
|\bmap(a)\cap X|= 0\ ,\\ \hat{0}_{\Ant(P)}, & \mbox{if
$|\bmap(a)\cap X|> 0$}\ ,
\end{array}\right. \\
\{a\}/ X&:= \left\{\begin{array}{ll} \{a\}, & \mbox{if }
|\bmap(a)\cap X|= 0\ ,\\ \bmap(\bmap(a)-X), & \mbox{if
$|\bmap(a)\cap X|> 0$ and $\bmap(a)\not\subseteq X$}\ ,\\
\hat{1}_{\Ant(P)}, & \mbox{if $\bmap(a)\subseteq X$}\ .
\end{array}\right.
\end{align*}

\item
If $A$ is a nontrivial antichain in $P$ then the order-theoretic
{\em deletion\/} $A\backslash X$ and {\em contraction\/} $A/ X$ of
$A$ in $P$ are the following joins in $\Ant(P)$:
\begin{equation}
\label{Label15} A\backslash X :=\bigvee_{a\in A} (\{a\}\backslash
X )\ ,\ \ \ A/ X :=\bigvee_{a\in A} (\{a\}/ X)\ .
\end{equation}

\item
The order-theoretic {\em deletion\/} and {\em contraction\/} of
the trivial antichains in $P$ are:
\begin{equation*}
\hat{0}_{\Ant(P)}\backslash X=\hat{0}_{\Ant(P)}/
X:=\hat{0}_{\Ant(P)}\ ,\ \ \ \hat{1}_{\Ant(P)}\backslash
X=\hat{1}_{\Ant(P)}/ X:=\hat{1}_{\Ant(P)}\ .
\end{equation*}
\end{itemize}

The map $(\backslash X):\Ant(P)\to\Ant(P)$, $A\mapsto A\backslash
X$, is called the {\em operator of deletion\/} on $\Ant(P)$; it is
a coclosure operator on $\Ant(P)$. The map $(/
X):\Ant(P)\to\Ant(P)$, $A\mapsto A/ X$, is called the {\em
operator of contraction} on $\Ant(P)$; it is a closure operator on
$\Ant(P)$~\cite[Theorem~2.5]{M2}.

Let the poset $L$ from Theorem~\ref{Label7} be the lattice
$\Ant(P)$. In this context, the maps $\bmap, (\backslash X), (/X)
:\Ant(P)\to\Ant(P)$ are instances of the maps $\beta$, $\delta$,
and $\gamma$ from Theorem~\ref{Label7}, respectively. In
particular, (\ref{Label9}) and~(\ref{Label11}) read as follows:

\begin{lem}
\label{Label17} For any antichain $A$ in $P$, the relations
\begin{align}
\label{Label13} \bmap(\bmap(A)\backslash X)&\geq A/ X\ ,\\
\label{Label14} \bmap(\bmap(A)/ X)&\geq A\backslash X
\end{align}
hold in $\Ant(P)$.
\end{lem}

\begin{proof}
There is nothing to prove if $A$ is trivial.

Let $\{a'\}$ be a nontrivial one-element antichain in $P$.
\begin{itemize}
\item[1.] Suppose that $|\bmap(a')\cap X|=0$. In this case we have
\begin{equation*}
\bmap\bigl(\bmap(a')\backslash
X\bigr)=\bmap\bigl(\bmap(a')\bigr)\geq \{a'\}=\{a'\}/ X
\end{equation*}
and
\begin{equation*}
\bmap\bigl(\bmap(a')/ X\bigr)=\bmap\bigl(\bmap(a')\bigr)\geq
\{a'\}= \{a'\}\backslash X\ .
\end{equation*}
\item[2.] Suppose that $|\bmap(a')\cap X|> 0$ and $\bmap(a')\not\subseteq X$.
In this case we have
\begin{equation*}
\bmap\bigl(\bmap(a')\backslash X\bigr)=\bmap\bigl(\bmap(a')-
X\bigr)=\{a'\}/ X
\end{equation*}
and
\begin{equation*}
\bmap\bigl(\bmap(a')/
X\bigr)=\bmap\bigl(\hat{1}_{\Ant(P)}\bigr)=\hat{0}_{\Ant(P)}=
\{a'\}\backslash X\ .
\end{equation*}
\item[3.] If $\bmap(a')\subseteq X$ then we have
\begin{equation*}
\bmap\bigl(\bmap(a')\backslash X\bigr)=\bmap\bigl(\bmap(a') -
X\bigr)=\bmap(\hat{0}_{\Ant(P)})= \hat{1}_{\Ant(P)}=\{a'\}/ X
\end{equation*}
and
\begin{equation*}
\bmap\bigl(\bmap(a')/
X\bigr)=\bmap\bigl(\hat{1}_{\Ant(P)}\bigr)=\hat{0}_{\Ant(P)}=\{a'\}\backslash
X\ .
\end{equation*}
\end{itemize}

Now, let $A$ be an arbitrary nontrivial antichain in $P$. On the
one hand, we by definition~(\ref{Label16}) have
\begin{equation*}
\bmap\bigl(\bmap(A)\backslash
X\bigr)=\bmap\Bigl(\left(\bigwedge_{a\in
A}\bmap(a)\right)\backslash X\Bigr)\ \ \text{and}\ \
\bmap\bigl(\bmap(A)/ X\bigr)=\bmap\Bigl(\left(\bigwedge_{a\in
A}\bmap(a)\right)/ X\Bigr)
\end{equation*}
in $\Ant(P)$. On the other hand, for any element $a'\in A$, we
have
\begin{equation*}
\bmap\Bigl(\left(\bigwedge_{a\in A}\bmap(a)\right)\backslash
X\Bigr)\geq\bmap\bigl(\bmap(a')\backslash X\bigr)\geq\{a'\}/ X
\end{equation*}
and
\begin{equation*}
\bmap\Bigl(\left(\bigwedge_{a\in A}\bmap(a)\right)/
X\Bigr)\geq\bmap\bigl(\bmap(a')/ X\bigr)\geq\{a'\}\backslash X
\end{equation*}
in $\Ant(P)$. Definitions~(\ref{Label15}) now imply
relations~(\ref{Label13}) and~(\ref{Label14}).
\end{proof}

With the help of relation~(\ref{Label8}) and Lemma~\ref{Label17},
we come to the following conclusion:
\begin{cor}[\cite{M2}, Theorem~2.6]
For any antichain $A$ in $P$, the relation
\begin{equation*} \bmap(A)\backslash X\ \leq\
\bmap(A/X)\ \leq\ \bmap(A)\ \leq\ \bmap(A)/X\ \leq\
\bmap(A\backslash X)
\end{equation*}
holds in $\Ant(P)$.
\end{cor}



\begin{thebibliography}{00}
\bibitem{Aigner}
M.~Aigner, {\em Combinatorial Theory}, Grundlehren der
Mathematischen Wissenschaften, vol.~234, Springer (1979).

\bibitem{BjH}
A.~Bj\"{o}rner and A.~Hultman, {\em A note on blockers in posets},
Ann.~Comb., {\bf 8} (2004), 123--131.

\bibitem{BPS}
A.~Bj\"{o}rner, I.~Peeva and J.~Sidman, {\em Subspace arrangements
defined by products of linear forms}, J.~London~Math.~Soc.~(2),
{\em to appear}.

\bibitem{CFM}
R.~Cordovil, K.~Fukuda and M.L.~Moreira, {\em Clutters and
matroids}, Discrete~Math., {\bf 89} (1991) 161--171.

\bibitem{C}
G.~Cornu\'{e}jols, {\em Combinatorial Optimization. Packing and
Covering}, CBMS-NSF Regional Conference Series in Applied
Mathematics, vol.~74, Philadelphia PA: SIAM (2001).

\bibitem{CH}
Y.~Crama and P.L.~Hammer, with contributions by C.~Benzaken,
J.C.~Bioch, E.~Boros, N.~Brauner, V.~Gurvich, L.~Hellerstein,
T.~Ibaraki, A.~Kogan, K.~Makino, R.~P\"{o}schel, I.~Rosenberg,
B.~Simeone and B.~Vettier, {\em Boolean Functions}, in
preparation.

\bibitem{EF}
J.~Edmonds and D.R.~Fulkerson, {\em Bottleneck extrema},
J.~Combinatorial~Theory, {\bf 8} (1970) 299--306.

\bibitem{GLS}
M.~Gr\"{o}tschel, L.~Lov\'{a}sz and A.~Schrijver, {\em Geometric
Algorithms and Combinatorial Optimization} Algorithms and
Combinatorics, vol.~2, Springer (1993).

\bibitem{L}
A.~Lehman, {\em A solution of the Shannon switching game},
J.~Soc.~Indust.~Appl.~Math., {\bf 12} (1964) 687--725.

\bibitem{M1}
A.O.~Matveev, {\em On blockers in bounded posets},
Int.~J.~Math.~Math.~Sci., {\bf 26} (2001), no.~10, 581--588.

\bibitem{M2}
A.O.~Matveev, {\em A note on operators of deletion and contraction
for antichains}, Int.~J.~Math.~Math.~Sci., {\bf 31} (2002),
no.~12, 725--730.

\bibitem{M3}
A.O.~Matveev, {\em Extended blocker, deletion, and contraction
maps on antichains}, Int.~J.~Math.~Math.~Sci., {\bf 2003} (2003),
no.~10, 607--616.

\bibitem{PDS}
P.D.~Seymour, {\em The forbidden minors of binary clutters},
J.~London~Math.~Soc. (2) {\bf 12} (1975/1976) 356--360.

\bibitem{St}
R.P.~Stanley, {\em Enumerative Combinatorics, vol.~1}, Cambridge
Studies in Advanced Mathematics, vol~49, Cambridge University
Press (1997).
\end{thebibliography}
\end{document}